\documentclass[11pt]{amsart}
\usepackage{amssymb}
\usepackage{amsmath}
\usepackage{graphicx}
\usepackage{amsthm}

\newtheorem{lemma}{Lemma}

\theoremstyle{remark}

\numberwithin{equation}{section}

\newcommand{\calG}{\mathcal{G}}

\newcommand{\calL}{\mathcal{L}}

\newcommand{\calO}{\mathcal{O}}
\newcommand{\calP}{\mathcal{P}}
\newcommand{\calQ}{\mathcal{Q}}

\newcommand{\calX}{\mathcal{X}}

\newcommand{\calI}{\mathcal{I}}

\newcommand{\bbF}{\mathbb{F}}
\newcommand{\bbC}{\mathbb{C}}

\newcommand{\bbP}{\mathbb{P}}

\newcommand{\bbZ}{\mathbb{Z}}
\newcommand{\bbG}{\mathbb{G}}

\def\Aut{\text{Aut}}

\def\Pic{\text{Pic}}

\def\Bir{\text{Bir}}
\def\El{\text{El}}
\def\cre{\text{cr}}
\def\Ell{\text{El}(n;p_1,\ldots,p_m)}
\def\ass{\text{ass}}
\address{Department of Mathematics, University of Michigan, Ann
Arbor, MI 48109, USA}
\email{idolga@umich.edu}
\thanks{Research partially supported by NSF grant DMS 990780}
\dedicatory{To the memory of Fabio Bardelli}
 
\title{On certain families of elliptic curves in projective space}
\author{Igor V. Dolgachev}
\begin{document}
\maketitle

%\begin{abstract}
%\end{abstract}

\section{Introduction} Let $\El(n;p_1,\ldots,p_m)$ be the family of elliptic 
curves of degree 
$n+1$ in $\bbP^n$ containing a fixed set of $m$ distinct points 
$p_1,\ldots,p_m$. In  modern language, 
$\El(n;p_1,\ldots,p_m)$ is  a  Zariski-open subset of a fibre of the 
evaluation map $\text{ev}:M_{1,m}(\bbP^n,n+1)\to (\bbP^n)^m$. Its general 
fibre, if not empty, is  an irreducible variety of dimension $(n+1)^2-m(n-1)$. The largest possible $m$ for 
which  $\text{ev}$ is dominant is equal to $9 (n= 2,5), 8 (n= 3, 4),  n+3 (n\ge 6)$. In the last case we 
show that the general  fibre is isomorphic to an open subset of a complete intersection of $n-2$ diagonal 
quadrics in $\bbP^{n+2}$. In particular, birationally, it is a Fano variety if $n\le 6$, a Calabi-Yau if 
$n = 7$, and of general type if $n\ge 8$. The group $\calG_n = (\bbZ/2\bbZ)^{n+2}$ acts naturally in 
$\bbP^{n+2}$ by multiplying the projective coordinates with $\pm 1$. The corresponding action of a subgroup 
of index 2 of $\calG_n$  is induced by a certain group of Cremona transformations in 
$\bbP^n$ which we explicitly describe. 

There are three cases when $\Ell$ is of expected dimension 0: $(n,m) = (2,9),(3,8), (5,9)$. It is well known that in the first two cases $\Ell$ consists of one point. Less known is the fact that the same is true in the case $(5,9)$. D. Babbage \cite{Babbage} attributes this result to T. G. Room. Apparently it was proven much earlier by A. Coble \cite{CobleAss}. We reproduce Coble's proof in the paper. This result implies the existence of a rational elliptic fibration $f:\bbP^5 - \to \El(5,p_1,\ldots,p_8)$ which is an analog of the well-known rational elliptic fibrations 
$\bbP^2 - \to \bbP^1 = \El(2;p_1,\ldots,p_8)$ and $\bbP^3 - \to \bbP^2 = \El(3;p_1,\ldots,p_7)$. We show that its locus of points of indeterminacy is a certain 3-fold, a  Weddle variety  studied intensively by Coble \cite{CobleWeddle},
\cite{CobleChantler}, \cite{CobleBul}.

I am grateful to D. Eisenbud, S. Mukai and R. Vakil for valuable discussions on the topic of this paper.

\section{Association}\label{s1} 
\subsection{}\label{association} We start with reminding the classical theory of association of finite sets of points. We follow the modern exposition of this theory given in \cite{Eisenbud}. 

Let $Z$ be a Gorenstein scheme of dimension 0 over a field $k$,  $\calL$ be an invertible sheaf on $Z$ and $V\subset H^0(Z,\calL)$ be a linear system. The duality pairing
\[H^0(Z,\calL)\times H^0(Z,\omega_Z\otimes \calL^{-1})\to H^0(Z,\omega_Z)\xrightarrow{\text{trace}} k\]
allows one to define the subspace $V^\perp\subset H^0(Z,\omega_Z\otimes \calL^{-1})$. The pair \[(V^\perp,\omega_Z\otimes \calL^{-1})\]
 is called the \emph{Gale transform} of $(V,\calL)$. 

\subsection{}\label{new} Assume that $Z$ is realized as a closed subscheme of a connected smooth $d$-dimensional scheme $B \subset \bbP^n$. Let $\calI_Z$ be its sheaf of ideals in $B$. We have a short exact sequence of sheaves on $B$
\[0\to \calI_Z(1) \to \calO_B(1)\to \calO_Z(1)\to 0,\]
which gives an exact sequence
\[H^0(B,\calO_B(1))\to H^0(Z,\calO_Z(1)) \to H^1(B,\calI_Z(1))\to H^1(B,\calO_B(1)).\]
Let us assume that $B$ is embedded by the complete linear system, i.e., the restriction map $H^0(\bbP^n,\calO_{\bbP^n})\to H^0(B,\calO_B(1))$ is bijective. Then the image of the map 
$H^0(B,\calO_B(1))\to H^0(Z,\calO_Z(1))$ is equal to $V$. If we also assume that $H^1(B,\calO_B(1)) = 0$, then we will be able to identify $V^\perp$ with $H^1(B,\calI_Z(1))^*$. 
Using Serre's duality, we get
\begin{equation}
\label{eq21 }
V^\perp \cong \text{Ext}^{d-1}(\calI_Z(1),\omega_B).\end{equation}
\subsection{}\label{ss2}In particular, taking $B = \bbP^n$, we can fix a basis in $\text{Ext}^{n-1}(\calI_Z(1),\omega_{\bbP^n})$, so that the linear system $\left |V^\perp\right |$ maps $Z$ to a subset $Z^{\ass}$ in the projective space 
\begin{equation}
\label{eq22 }
\bbP((V^\perp)^*) \cong \bbP(H^1(\bbP^n,\calI_Z(1))) \cong \bbP^{m-n-2}.
\end{equation}
The set  is called the \emph{associated set} of $Z$. It depends on a choice of a basis in  $H^1(\bbP^n,\calI_Z(1))$. If  $Z$ is a reduced scheme defined by an ordered set of $m$ distinct points 
$(p_1,\ldots,p_m)$ in $\bbP^n$, then the associated set $Z^{\ass}$ is an ordered set of points $(q_1,\ldots,q_m)$ in $\bbP^{m-n-2}$.

\subsection{}\label{ss23} Now let $B$ be a curve. Then $Z\subset B$ can be identified with a positive divisor, and 
\begin{equation}
\label{eq23 }
V^\perp \cong \text{Ext}^{d-1}(\calI_Z(1),\omega_B) = H^0(B,\calO_B(K_B+Z-H)),
\end{equation}
where $H$ is a hyperplane section and $K_B$ is a canonical divisor.
Thus we see that the associated scheme $Z^{\ass}$ is equal to the image of the divisor $Z$ under the map $B\to \bbP^{m-n-2}$ given by the complete linear system $\left |K_B+Z-H\right |$.

\subsection{Remark}\label{tyu} Let $Z= \{p_1,\ldots,p_m\}\subset \bbP^n$
considered as a reduced scheme. A choice of an order defines a basis in the 
space 
$H^0(Z,\calO_Z(1))$. Choosing  a basis in $H^0(\bbP^n,\calO_{\bbP^n}(1))$ we can represent $Z$ by the matrix 
$A$ of the linear map $H^0(Z,\calO_Z(1))^*\to H^0(\bbP^n,\calO_{\bbP^n}(1))^*$ which is equal to the 
transpose of the restriction map $H^0(\bbP^n,\calO_{\bbP^n}(1))\to H^0(Z,\calO_Z(1))$. The columns of 
$A$ are projective coordinates of the points $p_i$'s. Let $B$ be a matrix whose rows form a basis in the nullspace of $A$. Then the columns of $B$ can be chosen as the projective coordinates of the ordered set of points $q_1,\ldots,q_m$ in $\bbP^{m-n-2}$. These points represent the associated set of points $Z^{\ass}$ which is well defined up to projective equivalence. This is the original classical definition of association (see \cite{CobleAss}, \cite{Coble}). Let $P_n^m$ be the space of $m$
ordered points in
$\bbP^n$ modulo projective equivalence. This is a GIT-quotient with
respect to the democratic linearization (i.e the unique linearization
compatible with the action of the symmetric group $S_m$). The association defines an isomorphism of algebraic varieties (see \cite{DO}):
\[\ass:P_n^m\cong P_{m-n-2}^m.\]
As far as I know the first  cohomological interpretation of the association was given by A. Tyurin \cite{Tyurin}.

\section{Normal elliptic curves in $\bbP^n$ through $n+3$ points}\label{s3}
\subsection{}\label{ss31} We shall apply \ref{ss23} to the case when $B$ is an elliptic curve embedded in $\bbP^n$ by a complete linear system  (a normal elliptic curve of degree $n+1$). Let $\El(n;p_1,\ldots,p_m)$ be the set of normal elliptic curves in $\bbP^n$ containing $m\ge n+3$ distinct points $p_1,\ldots,p_m$. As in the Introduction we consider  $\El(n;p_1,\ldots,p_m)$ as an open Zariski subset of the fibre of the evaluation map $\text{ev}:M_{1,m}(\bbP^n,n+1)\to (\bbP^n)^m$ over $(p_1,\ldots,p_m)\in (\bbP^n)^m$. Let $Z$ be the reduced scheme $\{p_1,\ldots,p_m\}$ and let $B\in \El(n;p_1,\ldots,p_m)$. We fix a basis in $\bbP(H^1(\bbP^n,\calI_Z(1)))$ to identify $Z^{\ass} = \{q_1,\ldots,q_m\}$ with an ordered subset of $\bbP^{m-n-2}$. It follows from \ref{ss23} that $Z^{\ass}$ lies on the image of $B$ under the map $\phi$ given by the complete linear system $|Z-H|$ of degree $m-n-1$. Recall  that a choice of a basis in $H^1(\bbP^n,\calI_Z(1))$ defines a choice of a basis in $H^0(B,\calO_B(Z-H))$;  so the image of each $B\in \Ell$ lies in the same space $\bbP^{m-n-2}$ and contains $Z^{\ass}$. If $m > n+3$, the map $\phi:B\to \bbP^{m-n-2}$ is an embedding.  Since the association is the duality, we obtain 

\subsection{Theorem} \label{iso}{\it Assume $m > n+3$. Let $Z = \{p_1,\ldots,p_m\}$ and $Z^{\ass} = \{q_1,\ldots,q_m\}$. The association defines an isomorphism of algebraic varieies}:
\[\El(n;p_1,\ldots,p_m) \cong \El(m-n-2;q_1,\ldots,q_m).\]
 
\subsection{Remark}\label{self} In the case $m = 2n+2$ one can speak about self-associated schemes $Z$, that is, Gorenstein schemes such that $Z^{ass} = Z$ after appropriate choice of a basis in 
$H^1(\bbP^n,\calI_Z(1))$. In this case $\left |H\right | = \left |Z-H\right |$, so that $Z\subset B$ is cut out in $B$ by a quadric. The association isomorphism from Theorem \ref{iso} is of course the identity.

\subsection{} Let us consider the exceptional case $m = n+3$. In this case the map $\phi:B\to \bbP^{1}$ is of degree 2. It maps $p_1,\ldots,p_{n+3}$ to the associated set of points $q_1,\ldots,q_{n+3}$  in $\bbP^1$. Conversely, fix a set $\calQ = \{q_1,\ldots,q_{n+3}\}$ of distinct points in $\bbP^1$ and let $Z_\calQ$ be the corresponding $0$-dimensional scheme on $\bbP^1$. Take  a degree 2 map $\phi:B\to \bbP^1$ and choose $n+3$ points $\calP = \{p_1,\ldots,p_{n+3}\} \subset  B$ such that $p_i \in \phi^{-1}(q_i)$. Let $Z_\calP$ be the correspionding $0$-dimensional subscheme of $B$. Since 
\[\phi_*(\phi^*(\calO_{\bbP^1}(-1))\otimes \calO_{Z_\calP}) = \calO_{Z_\calQ}(-1),\]
we can identify the corresponding vector spaces of sections. The associated 
set $Z_\calQ^{\ass}$ lies on the image of $\bbP^1$ under the map 
$\bbP^1\to \bbP^n$ given by the linear system 
$\left |\calO_{\bbP^1}(Z_\calQ)\otimes \calO_{\bbP^1}(-1)\right |$. The 
linear system 
\[\left |\calO_B(p_1+\ldots+p_{n+3})\otimes \phi^*(\calO_{\bbP^1}(-1))\right
|\] embeds  $B$ into the same $\bbP^n$ with $Z_\calP^{\ass} = Z_\calQ^{\ass} = Z$. Thus $(B,\calP)$ defines a point in $\El(n;Z)$.    The double cover $\phi:B\to \bbP^1$ is defined, up to isomorphism, by a choice of the branch divisor $W\subset \bbP^1$, i.e. an unordered set of 4 distinct points in $\bbP^1$; the latter can be identified with a point in $\bbP^4\setminus \Delta$, where $\Delta$ is the discriminant hypersurface. We will continue to identify $\bbP^k$ with the symmetric product $(\bbP^1)^{(k)}$. Let $\iota:B\to B$ be the covering involution. It is clear that the $(B,\calP)$ and $(B,\iota(\calP))$ define the same point in $\El(n;Z)$. Summarizing we obtain the following:

\subsection{Theorem}\label{iso2}{\it  Let $q_1,\ldots,q_{n+3}$ be distinct points in $\bbP^1$ and let $p_1,\ldots,p_{n+3}$ be the associated set of points in $\bbP^{n}$. The variety $\El(n;p_1,\ldots,p_{n+3})$ is isomorphic to the fibre of the evaluation map 
$\text{ev}: M_{1,n+3}(\bbP^1,2)\to (\bbP^1)^{n+3}$ over $(q_1,\ldots,q_{n+3})$. 
The forgetting map  $M_{1,n+3}(\bbP^1,2)\to M_{1}(\bbP^1,2) 
\cong \bbP^4\setminus \Delta$ restricted to the fibre is a Galois cover 
with the Galois group $\calG_n$ isomorphic to $(\bbZ/2\bbZ)^{n+2}$. }

\subsection{}  Fix a set $\calP$ of 
$m$ points $p_1,\ldots,p_m$ in $\bbP^1$ and consider the subvariety 
\[X = \{A\in \bbP^k: A \cap \calP \ne \emptyset\}.\]
It is clear that $X = X_1\cup\ldots \cup X_m$, where 
$X_i = \{A\in \bbP^k: p_i\in A\}.$ This is the image in $(\bbP^1)^{(k)}$ of a coordinate hyperplane in 
$(\bbP^1)^k$, and hence is a hyperplane in $\bbP^k$ of effective divisors $D$ of degree $m$  such that 
$D-p_i > 0$. Let us identify $\bbP^k$ with $|\calO_{\bbP^1}(k)|$ and let $R = v_k(\bbP^1)$ be the Veronese 
curve in $\left |\calO_{\bbP^1}(k)\right |^* = \check{\bbP}^k$.  Then $X_i$ is the hyperplane in $\left |\calO_{\bbP^1}(k)\right |$ corresponding to the point $v_k(p_i)\in R$. Thus $X$ is the union of hyperplanes corresponding to $m$ points on $R$. 

\subsection{} The previous discussion shows that the branch divisor of the map $\El(n;p_1,\ldots,p_{n+3})\to \bbP^4\setminus \Delta$ consists of the union of open subsets of $n+3$ hyperplanes corresponding to $n+3$ points on a rational normal curve in $\bbP^4$.  Since the fundamental group of the complement of $N$ hyperplanes in general linear position in $\bbP^k$ is isomorphic to $\bbZ^{N-1}$, the Galois $(\bbZ/2\bbZ)^{n+2}$-cover $X\to \bbP^4$ branched along the union of $n+3$ hyperplanes $H_i$ is defined uniquely up to isomorphism. Thus $\El(n;p_1,\ldots,p_{n+3})$ is an open subset of $X$. A well-known way to construct $X$ is as follows. Let 
\[l_i = \sum_{s=0}^4a_{is}t_i = 0,\quad i = 0,\ldots,n+2\]
 be linear equations defining the hyperplanes $H_i$. We assume that $l_i = t_i,  i = 0,\ldots,4$.
 Consider the map 
\[i:\bbP^4\to \bbP^{n+2}, (t_0,\ldots,t_4)\mapsto (l_0(t_0,\ldots,t_4),\ldots,l_{n+2}(t_0,\ldots,t_4)).\]
 Let $(x_0,\ldots,x_{n+2})$ be  projective coordinates in $\bbP^{n+2}$ and let $(y_0,\ldots,y_{n+2})$ be projective coordinates in another copy of $\bbP^{n+2}$. Consider the map $\psi:\bbP^{n+2}\to \bbP^{n+2}$ defined by the formula $(x_0,\ldots,x_{n+2}) = (y_0^2,\ldots,y_{n+2}^2)$. Then $X$ is isomorphic to the  pre-image of $i(\bbP^4)$ under the map $\psi$. It is easy to see that $X$ is a complete intersection of the quadrics
\[y_i^2-\sum_{s=0}^4a_{is}y_i^2 = 0,\quad  i = 5,\ldots, n+2.\]

\subsection{Corollary} $\El(n;p_1,\ldots,p_{n+3})$ is isomorphic to an open Zariski subset of a smooth complete intersection of $n-2$ quadrics in $\bbP^{n+2}$. In particular, $X$ is birationally of general type  for $n > 7$,  Calabi-Yau for $n = 7  $, and Fano for $n\le 6$.

\subsection{Remark} The projection from the last point defines a map 
\[\El(n;p_1,\dots,p_m)\to \El(n-1;p_1,\ldots,p_{m-1}).\]
 In the case $m = n+3$, this is a finite map of varieties of dimension 4. I do not know the degree of this map.

\section{Cremona action}\label{s4}
\subsection{}\label{ss41} Let $P_n^m$ be the space of $m$
ordered points in
$\bbP^n$ modulo projective equivalence (see Remark~\ref{tyu}).  Recall that
the Cremona action on $P_n^m$ is a homomorphism of groups:
\begin{equation}\label{eq41}
\cre_{n,m}:W_{n,m}\to \Bir(P_n^m),\end{equation}
where $W_{n,m}$ is the Coxeter group corresponding to the Coxeter diagram
defined by a $T_{2,n+1,m-n-1}$-graph. Let
$s_0,\ldots,s_{m-1}$ be its Coxeter generators with $s_0$ corresponding
to the vertex of the Coxeter diagram with the arm of length $2$.

\begin{center}
  \includegraphics[width=2in]{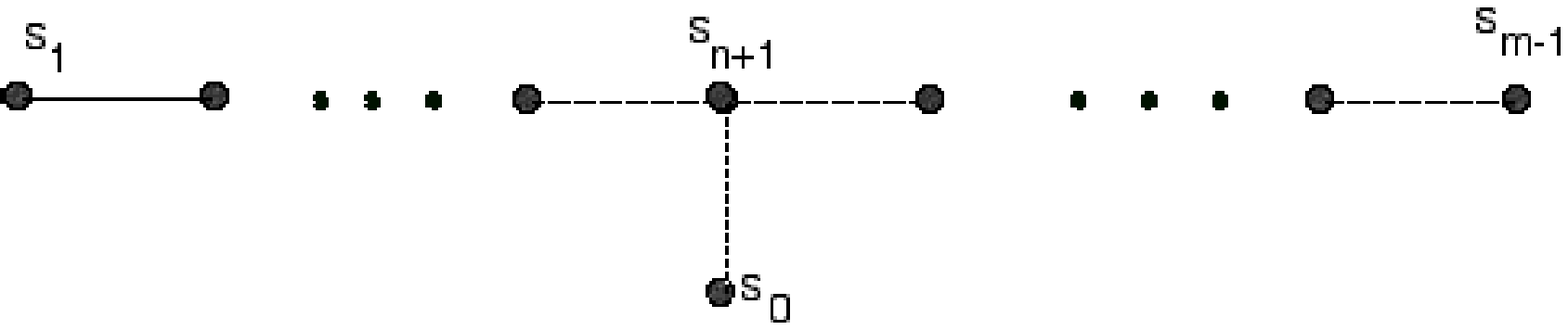}
\end{center}

 In the action \eqref{eq41}, the
element $s_0$ acts by means of the standard Cremona transformation 
\[T: \bbP^n - \to \bbP^n, \quad (z_0,\ldots,z_n)\mapsto (z_1\cdots z_n, z_0z_2\cdots z_n,\ldots, z_0z_1\cdots z_{n-1})\] 
 as follows. Choose an open subset $U$ of $P_n^m$ representing  stable orbits of sets of distinct points $\calP = (p_1,\ldots,p_m)$ such that $p_1,\ldots,p_{n+1}$ span $\bbP^n$. Then one can represent any point $x$ of $U$ with a  set $\calP =
(p_1,\ldots,p_m)$ such that 
\begin{equation}\label{eq42}
p_1 = (1,0,\ldots,0),\ldots, p_{n+1} =
(0,\ldots,0,1).\end{equation}
Then we define
\[\cre_{n,m}(s_0)(x) = (p_1,\ldots,p_{n+1},T(p_{n+2}),\ldots,T(p_m))\]
The subgroup $S_m$ of $W_{n,m}$  generated by $s_1,\ldots,s_{m-1}$ acts  via permutation of factors of $(\bbP^n)^m$. 

\subsection{}\label{ss42} For any subset $\calP\in (\bbP^n)^m$ representing a point $x\in U$  let $\pi_{\calP}:X(\calP)\to \bbP^n$ denote the blowing up of $\bbP^n$ with center at $\calP$. The Picard group $\Pic(X(\calP))$ has a natural basis $e_\calP= (e_0,e_1,\ldots,e_m)$, where
$e_0$ is the class of the pre-image of a hyperplane, and $e_i, i > 0,$ is the class of
the exceptional divisor blown-up from the point $p_i$. This basis is independent of a choice of a representative of the orbit $x$.  The group $W_{n,m}$ acts on
$\Pic(X(\calP))$ as follows. The subgroup $S_m$ acts by  permuting $e_i, i\ne 0,$ and
$s_0$ acts by:
\[e_0 \mapsto e_0-e_1-\ldots -e_{n+1}, \ e_i\mapsto e_0-e_1-\ldots -e_{n+1}+
e_i, i =
1,\ldots,n+1,\]
all $e_i$ are invariant for $i > n+1$.      
It is immediately checked that all elements $w\in W_{n,m}$ leave the anticanonical class
\[-K_{X(\calP)} = (n+1)e_0-(n-1)\sum_{i=1}^{n+3}e_i\]
invariant.
 
If $y =\cre_{n,m}(w)(x)$ and $\calQ$ is a representative of $y$, then there is a birational map 
$f: X(\calP)\to \calX(\calQ)$ which is an isomorphism in codimension $\ge 2$ such that 
\[w(e_\calP) = f^*(e_\calQ).\]
In particular, if $y = x$, we can choose $\calQ = \calP$ and obtain a pseudo-automorphism $g$ of $X(\calP)$ (i.e. birational automorphism which is an isomorphism in codimension $\ge 2$) such that $w = g^*$.

\subsection{}\label{ss43} Recall from Remark \ref{tyu} that the association defines an isomorphism of algebraic varieties
\[\ass_{n,m}:P_n^m\to P_{m-n-2}^m.\]
It  commutes with the Cremona action in the following sense. There is a natural isomorphism of the 
Coxeter groups $\tau:W_{n,m}\to W_{n,m-n-2}$ defined by the unique isomorphism of the   Coxeter diagrams 
$T_{2,n+1,m-n-1}$ and $T_{2,m-n-1,n+1}$ which leaves the vertex corresponding to $s_0$ fixed. We have
\begin{equation}
\label{eq43}
\cre_{n,m}(w)(x) = \cre_{m-n-2,m}(\tau(w))(\ass_{n,m}(x)).\end{equation}

\subsection{}\label{ss44} We will be  interested in the special case
$m= n+3$. In this case $W_{n,m}\cong W(D_{n+3})$, the Weyl group of
the root system of type $D_{n+3}$. It is known that
$W(D_{n+3})\cong G_n\rtimes S_{n+3}$, where $G_n \cong (\bbZ/2)^{n+2}$ is
generated by the element $w_1 = s_0\circ s_{n+2}$ and their
conjugates.   Let
$x\in P_n^{n+3}$ be represented by a point set $\calP$ as in \eqref{eq42}.
We also assume that
$p_{n+2} = (1,\ldots,1), p_{n+3} = (a_0,\ldots,a_n)$ with $a_i \ne 0$. Then
\[\cre_{n,n+3}(w_1)(x) = (p_1,\ldots,p_{n+1},p_{n+3}',p_{n+2}),\]
where 
$p_{n+3}' = (1/a_0,\ldots,1/a_n).$ Consider the projective transformation
\[g:\bbP^n\to \bbP^n, \quad (t_0,\ldots,t_n)\mapsto (a_0t_0,a_1t_1,\ldots,a_nt_n).\]
Then
$w(\calP) = g(\calP)$, and hence $\cre_{n,n+3}(x) = x$. This shows that $w_1$ acts trivially on $P_n^{n+3}$. In particular, the Cremona action $\cre_{n,m}$ has the kernel which is a normal subgroup of $W_{n,n+3}$ containing $w_1$. It is easy to see that the kernel coincides with $G_n$. For example, this follows from \eqref{eq43} since $P_n^{n+3} \cong P_1^{n+3}$. Thus we obtain that $G_n$ acts by pseudo-automorphisms  of  the blow-up $X(\calP)$.

\subsection{}\label{ss45} Following \cite{CobleWeddle} we show in more details how the group $G_n$ acts on $X(\calP)$. Let $V$ be the vector space over $\bbF_2$
consisting of  subsets of $[n+3] = \{1,\ldots,n+3\}$ of even cardinality with the addition defined by $A+B = A\cup B\setminus (A\cap B)$. 
Let $A$ be the affine space over $V$ consisting of subsets of odd
cardinality.  For
any
$I\in A$ of cardinality
$2k+1$ consider the following  divisor class in $X = X(\calP)$:
\begin{equation}\label{dc}
D_I = ke_0-(k-1)\sum_{i\in I}e_i-k\sum_{i\in
\bar{I}}e_i,\end{equation}
where $\bar{I}$ denotes the complementary subset of $[n+3]$. 

\begin{lemma}\label{L1} There is a unique isomorphism $h:V \to G_n, J\mapsto w_J$ such
that $w_{\{n+2,n+3\}} = s_0\circ s_{n+2}$, and
\begin{equation}
\label{ wj}
 w_J(D_I) = D_{J+I}, \quad \forall J\in V, \ \forall I\in A
\end{equation}
\end{lemma}

\begin{proof} First we check this directly for $w_{\{n+2,n+3\}}$:
\[D_{I+\{n+2,n+3\}} = s_0\circ s_{n+2}(D_I),\quad \forall I\in A.\]
Then we extend this to all subsets of cardinality 2 of $V$ by defining, for any $\sigma\in S_{n+3}$,
$w_{\sigma(\{n+2,n+3\})} = \sigma\circ w_{n+2,n+3}\circ \sigma^{-1}$ and using that
\[D_{I+\sigma(\{n+2,n+3\})} = D_{\sigma(\sigma^{-1}(I)+\{n+2,n+3\})} = \sigma\circ w_{n+2,n+3}\circ \sigma^{-1}(D_I).\]
 Finally, we extend the map $h$ to all even subsets by linearity. 
\end{proof}

\subsection{}\label{ss46} Let $J\in V$ with $\#J = 2k$. Then
\[w_J(e_s) = D_{J\cup \{s\}} = ke_0-(k-1)\sum_{i\in J}e_i-k\sum_{i\not\in J}e_i+ e_s,\]
 if $s\not\in J$, and 
\[w_J(E_s) = D_{J\cup \{s\}} = (k-1)e_0-(k-2)\sum_{i\in J}e_i-(k-1)\sum_{i\not\in J}e_i- e_s\]
otherwise. We have
\[-K_{X(\calP)} = (n+1)e_0-(n-1)\sum_{i=1}^{n+3}e_i = w_J((n+1)e_0-(n-1)\sum_{i=1}^{n+3}w_J(E_i)\]
\[= (n+1)w_J(e_0)-(n-1)(n+1)ke_0-(n-1)(k-1)(n+1)-1)\sum_{i\in J}e_i\]
\[-(n-1)(k(n+1)-1)\sum_{i\not\in J}e_i.\]
This gives
\begin{equation}
\label{eq48}
w_J(e_0) = (k(n-1)+1)e_0-(k-1)(n-1)\sum_{i\in J}e_i-k(n-1)\sum_{i\not\in J}e_i.\end{equation}
This shows that the pseudo-automorphism of $X_{\calP}$ induced by $\cre_{n,m}(w_J)$ is given by the linear system of hypersurfaces of degree $kn-k+1$ passing through the points $p_i\not\in J$ with multiplicity $(k-1)(n-1)$ and through the points $p_i\in J$ with multiplicity $k(n-1)$.

\subsection{}\label{ss47} Assume that $n = 2g-1$ is odd. Then we may take $J = [n+3]$ so that
\begin{eqnarray}
w_{[n+3]}(e_0) &= &(2g^2-1)e_0-2g(g-1)\sum_{i=1}^{n+3}e_i,\\
w_{[n+3]}(e_j) & = & gH-(g-1)\sum_{i=1}^{n+3}e_i-e_j 
\end{eqnarray}
\[= -\frac{1}{2}K_{X(\calP)}-e_j,\quad  j = 1,\ldots,n+3.\] 
It follows from Lemma~\ref{L1} that, for any $I\in A$,
\[w_{[n+3]}(D_I) = D_{\bar{I}},\]
and hence
\[w_{[n+3]}(D_I) + D_{\bar{I}} = ge_0-(g-1)\sum_{i=1}^{n+3}e_i\in \left |-\frac{1}{2}K_{X(\calP)}\right |.\]
The subgroup $G_n\subset W_{n,n+3}$ acts on $X_{\calP}$ and leaves the half-anticanonical linear system 
\[\left |-\frac{1}{2}K_{X(\calP)}\right | = \left |ge_0-(g-1)\sum_{i=1}^{n+3}e_i\right |\]
 invariant. The distinguished element $w_{[n+3]}\in G_n$ leaves  $2^{g}$ divisors $D_I+D_{\bar{I}}\in \left |-\frac{1}{2}K_{X(\calP)}\right | $ invariant.
Coble shows (\cite{CobleWeddle}) that 
\begin{equation}
\label{can}
\dim \left |-\frac{1}{2}K_{X(\calP)}\right |  = 2^g-1\end{equation}
and the divisors $D_I+D_{\bar{I}}$ span $\left |-\frac{1}{2}K_{X(\calP)}\right |$. In particular, $w_{[n+3]}$ 
acts identically in the projective space $\left |-\frac{1}{2}K_{X(\calP)}\right | $.

\subsection{}\label{ss48} Let $\tilde{W}_g\subset X_{\calP}$ be the locus of fixed points of $w_{[n+3]}$. This is defined as the closure of the fixed locus of $w_{[n+3]}$ restricted to an invaraint open subset where $w_{[n+3]}$. The projection $W_g$ of $\tilde{W}_g$ in $\bbP^n$ is a Weddle variety of dimension $g$. Coble proves that the linear system $\left |-\frac{1}{2}K_{X(\calP)}\right |$ maps $W_g$ to $\bbP^{2^g-1}$ and the image is isomorphic to the Kummer variety of the Jacobian variety $\text{Jac}(C)$ of  a hyperelliptic curve $C$ of genus $g$ embedded via $\left |2\Theta\right |$ on  $\text{Jac}(C)$.  The curve $C$ is the hyperelliptic curve corresponding to $n+3 = 2g+2$ points on $\bbP^1$ associated to $(p_1,\ldots,p_{n+3})$. 

\subsection{}\label{ss49} Let $p\in W_g$ be a fixed point of $w_{[n+3]}$. Consider the set of 
points $(p_1,\ldots,p_{n+3},p)$ and let  
$\calQ = (q_1,\ldots,q_{n+3},q)$ be the associated set in $\bbP^2$. We can view $w_{[n+3]}$ as an element of $W_{n,n+4}$. The projection from $q$ defines $n+3$ points in $\bbP^1$ which are associated to $(p_1,\ldots,p_{n+3})$ (see \cite{DO}, \cite{CobleAss}). Using \eqref{eq43}, we see that $\cre_{n,n+4}(w_{[n+3]})$ corresponds to $\cre_{2,n+4}(w_0)$ for some $w_0\in W_{2,n+4}$. Moreover, we may assume that $\cre_{2,n+4}(w_0)$ defines an automorphism $g$ (in dimension 2 any pseudo-automorphism is an automorphism) of the blow-up $X(\calQ)$. Since $\cre_{n,n+3}(w_{[n+3]})$ fixes the orbit of $(p_1,\ldots,p_{n+3})$, we see that $\cre_{2,n+4}(w_0)$ leaves the lines $\langle q_i\rangle$ invariant. This implies that the automorphism $g$ is induced by a de Jonqui\`eres transformation of the projective plane (see \cite{SR}, p.150). Its fixed set of points is a hyperelliptic curve of degree $g+2$ with $g$-multiple points at $q$ and tangent lines $\langle q,q_i\rangle$. This curve is isomorphic to the curve $C$ from \ref{ss48}. When we vary the point $p$ in the Weddle variety $W_g$ the subvariety of $P_2^{n+4}$ of associated point sets $\calQ$ consists of orbits of points sets $(q_1,\ldots,q_{n+3},q)$ such that $C$ admits a plane model $C'$ of degree $g+2$ with Weierstrass points at $q_1,\ldots, q_{n+3}$ and a $g$-multiple point at $q$. The $g$ branches at $q$ define an effective divisor $D$ of degree $g$ on $C$ such that the map $C\to C'$ is given by the linear system $|g_2^1+D|$, where $g_2^1$ is the linear series of degree 2 defined by the hyperelliptic involution.  In this way the Weddle variety $W_g$ becomes birationally isomorphic to the symmetric product $C^{(g)}$ modulo the hyperelliptic involution, and hence to the Kummer variety of $\text{Jac}(C)$. All of this can be found in \cite{Coble}, \S 38.

\subsection{}\label{ss410} Let us see how the group $G_n$, in its realization as the group of pseudo-automorphisms of the blow-up $X(\calP)$, acts on the variety $\El(n;p_1,\ldots,p_n)$. For any $B\in \El(n;p_1,\ldots,p_n)$ let $\bar{B}$ be  the proper inverse transform under the projection $\pi_\calP:X(\calP)\to \bbP^n$. Since $\calO_{X(\calP)}(e_0) \cong  \pi_\calP^*(\calO_{\bbP^n}(1))$, we see that $\bar{B}$ is embedded to $\bbP^n$ by the linear system $\left |e_0\right |$. It follows from \eqref{eq48} that 
\[w_J(\bar{B})\cdot e_0 =  \bar{B}\cdot w_J(e_0) \]
\[= \bar{B}\cdot [(k(n-1)+1)e_0-(k-1)(n-1)\sum_{i\in J}e_i-k(n-1)\sum_{i\not\in J}e_i]  \]
\[= k(n-1)+1)(n+1)-(k-1)(n-1)2k-k(n-1)(n+3-2k) = n+1.\]
This shows that  $B' = \pi_{\calP}(w_J(\bar{B}))$ is  
embedded in $\bbP^n$ by a complete 
linear system. Similarly we check that 
\[w_J(\bar{B})\cdot E_i = \bar{B}\cdot w_J(E_i) = 1,\]
and hence $B'$ passes through the points $p_1,\ldots,p_{n+3}$. Thus each $w_j$ acts on the variety 
$\El(n,p_1,\ldots,p_{n+3})$. We would like to verify that this action corresponds to the action of 
the Galois group $\calG_n$ of the cover $\El(n;p_1,\ldots,p_{n+3})\to \bbP^4$ from 
Theorem~\ref{iso2}. For this we consider the double covers $\phi:B\to \bbP^1$ and 
$\phi':B'\to \bbP^1$ defined by the linear systems 
$\left |\calO_B(-1)\otimes \calO_{B}(p_1+\ldots+p_{n+3})\right |$ and 
 $\left |\calO_{B'}(-1)\otimes \calO_{B'}(p_1+\ldots+p_{n+3})\right |$, respectively. Each maps the set $\calP = (p_1,\ldots,p_{n+3})$ to the same set of points in $\bbP^1$, the associated set of $\calP$. The isomorphism $w_J:B\to B'$ defines an isomorphism of the covers, and we have $\phi'(w_J(p_i)) = \phi'(p_i)$ for each $i = 1,\ldots,n+3$. Thus $B$ and $B'$ belong to the same fibre of the map 
$\El(n,p_1,\ldots,p_{n+3})\to \bbP^4$. 

The homomorphism  $G_n\to \calG_n$ which we have just constructed is not injective in the case when $n = 2g-1$ is odd. In fact, the linear system $\left |-\frac{1}{2}K_{X(\calP)}\right |$ cuts out a $g_2^1$ on each $\bar{B}$. Since $w_{[n+3]}$ acts identically on $\left |-\frac{1}{2}K_{X(\calP)}\right |$ it induces an automorphism of $\bar{B}$.Thus $B' = B$ in this case. It is clear that the fixed points of  $w_{[n+3]}$ on $\bar{B}$ are the points of intersection of $\bar{B}$ and the Weddle variety $W_g$. 
In particular, we obtain
\begin{equation}\label{eq49}
\# \bar{B}\cap W_g = 4
\end{equation}
unless all curves $b$ are conatined in $W_g$. However, by using the next Lemma and induction on $g$ we see that the latter does not happen.

\subsection{Lemma}\label{L2} {\it For any $J\in V$ with $\#J = 2k$, let $F_J$ be the set of fixed points of the involution $w_J$ in  $X_{\calP}$. Then $\dim F_J = k-1$. If $k > 2$, the projection  map $\pi_J$ from the linear space spanned by the points $p_i, i\in \bar{J}$ defines a finite map of degree $2^{n+2-2k}$ from $\pi_\calP(F_J)$ to the Weddle variety $W_{k-1}$ in $\bbP^{2k-3}$ corresponding to  the set $\pi_J(\{p_i\}_{i\in J})$. }

\begin{proof} See \cite{CobleWeddle},  p.456. 
\end{proof}

\subsection{Theorem} {\it The kernel of the homomorphism $G_n\to \calG_n$ is of order $\le 2$. It is  generated by $w_{[n+3]}$ if $n$ is odd.}

\begin{proof} Let $K$ be the kernel of the homomorphism $G_n\to \calG_n$ and $B\in\El(n;p_1,\ldots,p_{n+3})$. Assume $n$ is odd. Since $K$ leaves the linear system $\left |-\frac{1}{2}K_{X(\calP)}\right |$ invariant, $K$ leaves the $g_2^1$ on $B$ cut  by this linear system invariant. Thus the image of $K$ in $\Aut(B)$ is equal to the image of $w_{[n+3]}$. If $w_J\in K$ acts identically on $B$, then $B$ is contained in the fixed locus $F_J$ of $w_J$. Hence the 
union of all curves from $\El(n;p_1,\ldots,p_{n+3})$ is contained in $F_J$.  Let $\# J = 2k$. Since $\dim F_J$ is too small for $k \le 2$ in order this could happen, we may assume that $k > 3$. Projecting from the subspace spanned by the points $p_i, i\in \bar{J}$, we obtain that the Weddle variety $W_{k-1}$ contains all elliptic curves from $\El(2k-3;q_1,\ldots, q_{2k})$. But this contradicts \eqref{eq49}. 

Assume $n$ is even. Assume $B$ is preserved under some $w_J$. Since $w_j$ sends the set of points $p_1,\ldots,p_{n+3}$ to projectively equivalent set of points, we may assume that $w_J$ induces a projective transformation of $B$. Since $B$ is embedded in $\bbP^n$ by a linear system of odd degree, a projective automorphism of $B$ of order 2 cannot be a translation by a 2-torsion point (when we fix a group law on $B$). It must fix 4 points on $B$. Now we finish as in the previous case.
\end{proof}

\subsection{Remark} In the case $n$ is even, the homomorphism $G_n\to \calG_n$ must be injective. However I cannot prove it.

\section{A rational elliptic fibration on $\bbP^5$}

\subsection{} Now we consider the special case $n= 5, m= 9$. By Theorem~\ref{iso}
\[\El(5;p_1,\ldots,p_9) \cong \El(2;q_1,\ldots,q_9),\]
where $\calQ = (q_1,\ldots,q_9)$ is the point set associated to 
$\calP = (p_1,\ldots,p_9)$. It is well known that there is a unique nonsingular  cubic through a general set of 9 points in $\bbP^2$. Thus we obtain the following

\subsection{Theorem} (A. Coble)\label{coble} {\it Let $\{p_1,\ldots,p_9\}$ be a general set of 9 points in $\bbP^5$. Then there exists a unique elliptic curve of degree 6 containing these points. }

\subsection{} For a special set of 9 points $\calQ$ in $\bbP^2$, the dimension of the linear system of cubics through $\calQ$ varies from $0$ to $5$. However, if we assume that $\El(2;q_1,\ldots,q_9)\ne \emptyset$, then $\dim \El(2;q_1,\ldots,q_9) \le 1$.

\subsection{Proposition}\label{eis} {\it Assume that $\dim \El(2;q_1,\ldots,q_9) = 1$. Then the linear system defining the associated set of points  $\calP = (p_1,\ldots,p_9)$ in $\bbP^5$ is equal to the restriction of the complete linear system of conics. In particular, $\calP$ is equal to the image of $\calQ$ under a Veronese map $v_2:\bbP^2\to \bbP^2$.}

\begin{proof} The scheme $Z = \{q_1,\ldots,q_9\}$ is cut out by two cubics, and so is a an arithmetically Gorenstein, i.e., its homogeneous coordinate ring is Gorenstein. Now the assertion follows from Corollary 2.6 of \cite{Eisenbud}.
\end{proof}

\subsection{Remark} The assertion is claimed by Coble in \cite{CobleAss}. The proof uses a method which has not been justified so far (see Problem 2.3 in  \cite{Eisenbud}). 

\subsection{Lemma}(A. Coble)\label{ver} {\it Let $E$ be an elliptic curve of degree 6 in $\bbP^5$. Then $E$ is contained in exactly four Veronese surfaces.}

\begin{proof} Take 9 points on $E$ such that the associated set of points in $\bbP^2$ lies on a pencil of nonsingular cubics $C_t$. By Proposition~\ref{eis} and Theorem ~\ref{iso}, $E$ is the image of some $C_t$ under a Veronese map. Thus $E$ is contained in some Veronese surface (this can also be verified by counting constants). Let $V$ be a Veronese surface containing $E$. Then the linear system $\left |D\right |$ on $E$ cut out by the complete linear system $L$, where $L$ is an effective generator of $\Pic(V)$ satisfies $\left |2D\right | = \left |\calO_E(1)\right |$. Also, it comes with a fixed isomorphism $\phi:S^2(H^0(\calO_E(D)))\to H^0(\calO_E(2D))$. 
  Suppose two Veronese surfaces $V$ and $V'$ contain $E$ and cut out the same linear system $\left |D\right |$ on $E$. Then  the two embeddings $E\to V, E\to V'$ differ by a linear  endomorphism $T$ of $H^0(E,\calO_E(D))$ such that $S^2(T)$ is the identity. This obviously implies that $T$ is the identity, and hence $V = V'$.  On other hand, if $\left |D\right | \ne \left |D'\right |$, we can choose an isomorphism 
$S^2(H^0(E,\calO_E(D)))\to S^2(H^0(E,\calO_E(D')))$ such that the compositions of two embeddings $E\to |\calO_E(D)|^*$ and $E\to |\calO_E(D')|^*$ with the corresponding Veronese maps 
$ \left |\calO_E(D)\right |^*\to\left  |\calO_E(2D)\right |^*$ and $\left  |\calO_E(D')\right |^*\to \left |\calO_E(2D')\right |^*$ are the same. 
This proves the assertion. 
\end{proof}

\subsection{} Let $E\in \El(5,p_1,\ldots,p_8)$. Let $V$ be  a Veronese surface containing $E$. We may assume that $V$ is the 
image of the Veronese map $v_2:\bbP^2\to \bbP^5$ and $E = v_2(C)$ for some nonsingular plane cubic $C$. Let 
$q_1,\ldots,q_8$ be the points on $C$ such that $v_2(q_i) = p_i, i = 1,\ldots,8$. Choose the ninth point $q_9$ such that $\calQ = \{q_1,\ldots,q_9\}$ is the  base locus of a pencil of cubics. Then, by Proposition~\ref{eis} the associated set of points $\calQ^{\ass}$ is the image under a Veronese map. We can fix the Veronese map by requiring that $v_2(p_i) = q_i, \ i = 1,\ldots,8$. Let $v_2(q_9) = p_9$. Then by Theorem~\ref{iso}, $\El(5;p_1,\ldots,p_9) \cong \El(2;q_1,\ldots,q_9)$ is an open subset of $\bbP^1$. We want to show that the point $p_9$ is on the Weddle variety $W_3$ corresponding to the set of points $(p_1,\ldots,p_8)$. 

\subsection{Lemma}\label{L3} {\it Let $(q_1,\ldots,q_9)$ be an ordered set of base points of an irreducible pencil of cubic curves. Then there exists a unique plane curve $C$ of degree $5$ with a triple point at $q_9$ such that the lines $\langle q_i,q_9\rangle$ are tangents to $C$ at $q_i$. Conversely, given such a curve $C$, there exists a pencil of cubic curves with base points at $q_1,\ldots,q_9$.}

\begin{proof} We give a modern  version of the proof given by Coble in \cite{Coble}, \S 38.  Let $X$ be the blow-up of the points $q_1,\ldots,q_9$. This has a structure of an elliptic surface together with the zero section defined by the exceptional curve $E_9$ blown up from $q_9$. Let $F$ be the divisor class of a fibre and $S$ be the divisor class of the proper transform of a line through the point $q_9$. Then $\left |F +S\right |$ defines a double cover $f:X\to \bbP^1\times\bbP^1$ such that the pre-image of the first ruling is equal to the pencil $\left |F\right |$, and the pre-image of the second ruling is equal to the pencil $|S|$. The branch curve $B$ of $f$ is a curve of bi-degree $(4,2)$. It is a hyperelliptic curve of genus 3. Its pre-image in $X$ is equal to $2R$, where $R$ is isomorphic to $B$. On each nonsingular fibre   $F$ there are two ramification points of the pencil $\left |S\right |$ restricted to $F$. The curve $R$ is equal to the locus of  these points when $F$ varies in the pencil. Similarly, for each member of $\left |S\right |$ there are two ramification points of the restriction of the pencil $\left |F\right |$, and these ramification points are also on $R$. 

Let us show that the image of $R$ in $\bbP^2$ is the hyperelliptic curve $C$ from the assertion of the lemma. The image $\bar{E}_9$ of $E_9$ is a curve of bi-degree $(1,1)$ such that $f^{-1}(\bar{E}_9) = E_9+E_9'$, where $E_9'$ is the image of $E_9$ under the involution of $X$ defined by the double cover $f$. Since $\bar{E}_9$ is tangent everywhere to the branch curve $R$, we  see that $E_9$ intersects $R$ at 3 points. Thus $C$ has a triple point at $q_9$. Let $E_i$ be the exceptional curve blown up from $q_i$ and let $S_i$ be the proper inverse transform of the line  $\langle q_i,q_9\rangle $. Then $E_i+S_i\in \left |S\right |$, and $f(E_i) = f(S_i)$ is a fibre of the second ruling of $\bbP^1\times \bbP^1$ which is tangent to $B$ at a point $P$ such that $f^{-1}(P) = E_i\cap S_i$. This shows that $E_i\cap S_i\cap R\ne \emptyset$, and hence $C$ is tangent to $\langle q_i,q_9\rangle$ at $q_i$. Thus $C$ satisfies all the properties from the assertion of the lemma. Let us show its uniqueness. Suppose $C'$ is another curve satisfying the same properties. Let $R'$ be its proper inverse transform on $X$. We have $R\cdot R' = 25-9-16 = 0$, and $R^2 = R'{}^2 = 25 -9-8 = 8$. This contradicts the Hodge Index Theorem.

It remains to prove the converse. We will be brief letting the reader to fill
in the details. Let $R$ be the proper inverse transform of $C$ in $X$. Then the linear system $|R-S|$ defines a double cover $f:X\to \bbP^1\times \bbP^1\subset \bbP^3$. If the restriction of $f$ to $S$ is of degree 1, then its image is a plane section of $\bbP^1\times \bbP^1$ and hence   $|R-2S| = \{S'\}$ for some curve isomorphic to $S$. However, $(R-2S)^2 = 4, (R-2S)\cdot K_X =  0$ imply  that the arithmetical genus of $S'$ is equal to 1. This contradiction proves that the image of $S$ must be a fibre of a ruling on $\bbP^1\times \bbP^1$. This implies that $|R-2S|$ is a pencil of elliptic curves (the pre-image of another ruling). Its image in $\bbP^2$ is a pencil of cubics with base points $q_1,\ldots,q_9$.
\end{proof}

\subsection{Theorem}\label{last} {\it Let $p_1,\ldots,p_8$ be general points in $\bbP^5$. Then the locus of point $p_9\in \bbP^5$ such that $\dim \El(5;p_1,\ldots,p_9)  = 1$  is equal to the Weddle variety $W_3$ associated to the points  $p_1,\ldots,p_8$. }

\begin{proof} We know that if $\dim \El(5;p_1,\ldots,p_9)  = 1$, then $(p_1,\ldots,p_9)$ is associated to a base point set of a pencil of cubics in $\bbP^2$. Now the assertion follows from \ref{ss49} and   the previous Lemma.
\end{proof}

\subsection{Remark} The assertion of Theorem \ref{last} agrees  with \eqref{eq49} and Lemma \ref{ver}. We have four Veronese surfaces $V_i$ on $E\in \El(5;p_1,\ldots,p_8)$. In each $V_i$ we find the ninth point $p_9\in E$ such that there exists a pencil of cubics through $p_1,\ldots,p_9$. The  point $p_9$ must belong to $W_3$ and we have $E\cap W_3 = 4$.

\subsection{Remark} It was shown by Coble (\cite{CobleWeddle}, p.489) that $W_3$ is of degree 19 with singular points at $p_1,\ldots,p_8$ of multiplicity 9, triple lines $\langle p_i,p_j\rangle $ and the triple curve equal to the rational normal curve through the points $p_1,\ldots,p_8$.

\subsection{} By Theorem~\ref{coble} we have a rational elliptic fibration 
\[\Phi:\bbP^5 - \to \El(5;p_1,\ldots,p_8)\]
  defined by the universal family of elliptic curves passing through the points $p_1,\ldots,p_8$.  Its points of indeterminacy is the Weddle variety $W_3$. If we blow up $p_1,\ldots,p_9$ and then blow up the proper inverse transform $\tilde{W}_3$ of $W_3$ we expect to get a regular map $\tilde{\Phi}:\tilde{X} \to \El(5;p_1,\ldots,p_8)$. Let $\pi:\tilde{X}\to \bbP^5$ be the composition of these blow-ups. For each nonsingular points $p\in W_3\setminus \{p_1,\ldots,p_8\}$ the image of the line $\pi^{-1}(p)$ under $\tilde{\Phi}$  in $\El(5;p_1,\ldots,p_8)$ is the pencil of elliptic curves from $\El(5;p_1,\ldots,p_8)$ passing though $p$. Thus there are 4 rational curves through each general point in the variety $\El(5;p_1,\ldots,p_8)$. The restriction of $\tilde{\Phi}$ over each such curve is a rational elliptic surface corresponding to the associated set of points $\calQ$ of $\calP = (p_1,\ldots,p_8,p)$ in $\bbP^2.$ It is isomorphic to the proper inverse transform $\tilde{V}$ under $\pi$ of a Veronese surface containing $(p_1,\ldots,p_8,p)$. In the Cremona action $\cre_{2,9}:W_{2,9}\to \Bir(P_2^9)$ the orbit of $\calQ$ is fixed under the subgroup $G$ isomorphic to $H = \bbZ^8\rtimes (\bbZ/2\bbZ)$(see \cite{DO}, p. 124). Using \eqref{eq43}, we obtain that the orbit of the points set $\calP$ is fixed under the subgroup isomorphic to $H$ in the Cremona action $\cre_{5,9}$. In particular,  we find that the group of pseudo-automorphisms of the blow-up $X(\calP)$ contains an infinite group $G'$ isomorphic to $H$. Let $V'$ be the proper inverse transform of the Veronese surface $V$ in $X(\calP)$. The projection $\tilde{X}\to X(\calP)$ defines  an isomorphism  $\tilde{V} \to V'$.  The group $G$ leaves $V'$ invariant and its action on $V'$ corresponds  to the action of $G$ on $X(\calQ)$. Apparently (I have not checked it), the action of $
H$ on each $\tilde{V}$ is the restriction of the automorphism group of the generic fibre of the elliptic fibration $\tilde{\Phi}$ isomorphic to $H$. 

\subsection{Remark} The fact that the group of pseudo-automorphisms of the blow-up $X(\calP)$ is 
infinite  was used by S. Mukai (\cite{Mukai}) to construct an example of a linear action of the additive group $\bbG_a^3$ in the polynomial algebra $\bbC[T_1,\ldots,T_9]$ such that the  generated algebra of invariants is not  finitely generated.

%\bibliographystyle{plain}
	%\bibliography{mybibliography}

\end{document}